\documentclass[12pt]{amsart}
\usepackage{amssymb}

\def\Sp{\mathbb{S}}

\def\R{\mathbb{R}}

\def\tilde{\widetilde}
\def\epsilon{\varepsilon}
\def\phi{\varphi}
 \def\al{\alpha} \def\be{\beta} \def\ga{\gamma} \def\Ga{\Gamma}  
 \def\Om{\Omega}  \def\la{\lambda} 
 \def\De{\Delta}

\def\rank{{\rm rank}}

\def\pa{\partial }

\def\sm{\setminus }
\def\ss{\subset }
																																			
\def\beq{\begin{equation}}
\def\eeq{\end{equation}}

\newtheorem{theorem}{Theorem}[section]
\newtheorem{lemma}[theorem]{Lemma}
\newtheorem{proposition}[theorem]{Proposition}

\DeclareMathOperator{\Real}{Re}
\DeclareMathOperator{\Imaginary}{Im}

\begin{document}

\title[Minimal surfaces and boundary 
value problems]{Free boundary minimal surfaces and overdetermined boundary 
value problems}
\author{Nikolai  Nadirashvili}
\address[Nikolai S. Nadirashvili]{CNRS, I2M UMR 7353 --- Centre de Math\'ematiques et Informatique,
Marseille, France}
\email{nnicolas@yandex.ru} 
\author{Alexei V. Penskoi}\thanks{The work of the second author
was partially supported by the Simons Foundation and 
by the Young Russian Mathematics award.}
\address[Alexei V. Penskoi]{Department of Higher Geometry and Topology, 
Faculty of Mathematics and Mechanics, Moscow State University,
Leninskie Gory, GSP-1, 119991, Moscow, Russia \newline {\em and}
\newline Faculty of Mathematics,
National Research University Higher School of Economics,
6 Usacheva Str., 119048, Moscow, Russia \newline {\em and}
\newline Interdisciplinary Scientific Center
J.-V. Poncelet (ISCP, UMI 2615), Bolshoy Vlasyevskiy 
Pereulok 11, 119002, Moscow, Russia}
\email[Corresponding author]{penskoi@mccme.ru}
\date{}

\begin{abstract} In this paper we establish a connection between free boundary minimal surfaces in a ball
in $\mathbb{R}^3$ and free boundary  cones arising in a one-phase problem. 
We prove that a doubly connected minimal surface 
with free boundary in a ball is a catenoid.  
\end{abstract}
\maketitle

\section{Introduction}

In this paper we investigate free boundary minimal surfaces in a three-dimensional ball,
i.e. proper branched minimal immersions of a surface $M$ 
$$
u: M\longrightarrow\mathbb{B}^3\subset\mathbb{R}^3
$$
such that $u(M)$ meets the boundary sphere $\mathbb{S}^2=\partial\mathbb{B}^3$ orthogonally. It is a classical and developed subject, see e.g. the 
books~\cite{DierkesHildebrandtSauvigny2010,DierkesHildebrandtTromba2010a,DierkesHildebrandtTromba2010}. 
A celebrated result due to J.~C.~C. Nitsche~\cite{Nitsche1972} states that if $M$ is a disk then $u(M)$ is also
a plane disk.

Actually, in the paper~\cite{Nitsche1972}  a stronger result is announced. Namely, that this statement holds
for capillary surfaces and  the 
angle between $u(M)$ and $\mathbb{S}^2$ is a constant. Details of the proof
could be found in the paper~\cite{RosSouam1997}.

Recently, the result due to J.~C.~C. Nitsche was generalized by A.~Fraser and R.~Schoen 
in the paper~\cite{FraserSchoen2015} to the
case of a minimal disk satisfying the free boundary condition in a constant curvature ball of any dimension.

In this paper we prove that a free boundary doubly connected minimal surface  in a three-dimensional
euclidean ball is a piece of a catenoid. 

Since (i) a minimal map could be parametrized by a conformal
parameter, (ii) a double-connected domain is conformally equivalent to
an annulus
\begin{equation}\label{annulus}
A=\{z|\rho<|z|<1\}\subset\mathbb{C},
\end{equation}
and (iii) any ball could be transformed by a homotety to a unit ball,
it is sufficient to consider a map from annulus~\eqref{annulus} to the unit ball 
$\mathbb{B}^3_1\subset\mathbb{R}^3$
centered at the origin.

\begin{theorem}\label{th1}
Let  $u: \bar{A}\longrightarrow\mathbb{\bar{B}}^3_1\subset\mathbb{R}^3$
be a proper branched minimal immersion
such that $u\in C^2(A)\cap C^1(\bar{A}),$ $u(\partial A)\subset\mathbb{S}^2=\partial\mathbb{B}_1^3$
and $u(A)$ meets $\mathbb{S}^2$ orthogonally.
Then $u(A)$ is a part of a catenoid.
\end{theorem}  
This result was conjectured by A. Fraser and M. Li in~\cite{FraserLi2014}.

Recently A. Fraser and R. Schoen  proved the existence of free boundary minimal surfaces in 
$\mathbb{B}^3$ which have genus $0$ and $n$ boundary components, see the papers~\cite{FraserSchoen2016}, 
see also~\cite{FolhaPacardZolotareva2017} . 

Let us remark that A. Fraser and R. Schoen established in the papers~\cite{FraserSchoen2011,FraserSchoen2016} a remarkable connection between minimal surfaces with 
free boundaries in a ball and Riemannian metrics on surfaces with boundaries extremizing eigenvalues 
of the Steklov problem on these surfaces. Let us also remark that connections of spectral isoperimetry with 
minimal surfaces is known also for some other spectral problems, see
the paper~\cite{NadirashviliSire2015}.

%%
%Minimal surfaces with free boundaries is an important subclass of capillarily surfaces which are modeling physical problem of the behavior of an incompressible liquid in a container  in the absence of gravity, see the book of R. Finn    . The capillarily surfaces called stabled if it locally
%minimize a corresponding energy functional. A. Ros and R. Souam proved, ~\cite{RosSouam1997},  that
%the only capillarily stable and minimal surfaces in a ball of $\R^3$ are the totally geodesic disks or surfaces of genus $1$ with boundary having at most $3$ connected components. The existence of capillarily stable and minimal surfaces in a ball distinguished from the disk and the catenoid remains open.
%%

In this paper we  establish by means of the Minkowski transformation a connection between free boundary 
minimal surfaces and the extremal domains on the sphere $\Sp^2$ for the Dirichlet problem. The last spectral 
problem is related to the one-phase free boundary problem
for homogeneous functions defined in cones. By virtue of this connection we prove some new results for 
the one-phase free boundary problem.

The one-phase free boundary problem is a minimization of an integral
$$
J(v)=\int_{G\cap \{ v>0\} } (|\nabla v|^2+1)dx \rightarrow \min,
$$
where $v\geqslant 0$. It appears in some models related to the cavitational flow. If $v$ is a minimizer and 
$G^+ =\{ x\in G, v(x)>0\}$ then $v$ is a solution of the following overdetermined problem in $G,$
\begin{equation}\label{0} \left\{
\begin{array}{l l}
 v\geqslant  0&\mbox{in $G$}, \\
\De v=0 &\mbox{on  $G^+ $}, \\
|\nabla v|=C  &\mbox{on  $\pa G^+ \cap G$},\\
\end{array} \right. \end{equation}
where $C>O$ is a constant. H. W.  Alt and L. A. Caffarelli proved 
in the paper~\cite{AltCaffarelli1981} that the 
question of regularity of the boundary in the one-phase  free boundary problem 
could be reduced to the one-phase problem in a cone. Let $K\ss \R^n$ be an open 
($n$-dimensional) cone with a smooth lateral boundary. We are interested in the following overdetermined problem,

\begin{equation}\label{1} \left\{
\begin{array}{l l}
\De v=  0&\mbox{in $K$}, \\
v=0 &\mbox{on  $\pa K$}, \\
|\nabla v|=1  &\mbox{on  $\pa K\sm \{ 0\}$},\\
\end{array} \right. \end{equation}
where $v$ is a homogeneous degree $1$ function.  Let us  emphasise that the unknowns here
are both $v$ and $K.$ 
 One can define an energy for the solutions of  system~\eqref{1} related to $J$, 
see the papers~\cite{AltCaffarelli1981,CaffarelliSalsa2005}.
 A solution $v$ of  \eqref{1} is called stable if it is stable with respect to this energy. L.~A.~Caffarlelli,
 D.~Jerison and C.~E.~Kenig proved that in $\R^3$ the only stable solutions of \eqref{1} are linear functions, 
the correspondent cone $K$ is a half-space,
see the paper~\cite{CaffarelliJerisonKenig2004}. This result was extended to the dimension $4,$ see the paper~\cite{JerisonSavin2015}. On the other hand, D.~De Silva and D.~Jerison gave an example of a nontrivial energy minimizing solution in dimension $7$,
 see the paper~\cite{DeSilvaJerison2009}.
 
 We show that result from the paper~\cite{CaffarelliJerisonKenig2004}
holds if we just assume that
$K$ is a simply connected cone instead of the stability of the 
 solution of system~\eqref{1}.

 \begin{theorem}\label{th2} Suppose that $v$ and $K\subset\mathbb{R}^3$
is a solution of system \eqref{1}. Then
 
 a) if $K\cap \Sp^2$ is diffeomorphic to a disk then $K$ is a half-space;
 
 b) if $K\cap \Sp^2$ is diffeomorphic to an anulus then
 $\mathbb{R}^3\setminus K$
 is a circular cone formed by lines with aperture $2\arccos\tanh\al$, where $\al$ is a solution of the equation 
 $\sqrt\al\tanh\al=1$.
 \end{theorem}
 
 The proof of Theorem~\ref{th2} is based on the following involution on the space of homogeneous order $1$ functions. Let $f$ be a homogeneous function of
order $1$ defined in a cone $K\ss \R^3$. 
Consider the surface
 $H_f=\nabla f(K)$ called the h\'erisson of $f$ (the notion was introduced in
the paper~\cite{LangevinLevittRosenberg1988}). 
 The following theorem which goes back to Minkowski, 
see~\cite[$\S$ 78]{Blaschke1921}, holds:
 
 \begin{proposition}\label{prop} Let  $f$ be a homogeneous  function of degree $1$
 defined in a cone $K\ss \R^3$. Let $N(H_f)$ be a Gauss map of the h\'erisson $H_f$. Then 
 the map $N:H_f\longrightarrow\mathbb{S}^2$ is inverse at
 regular points of $H_f$ to 
 $\nabla f:\mathbb{S}^2\longrightarrow H_f$. Moreover at a regular point $z\in H_f$ the 
sum of curvature radii of the surface $H_f$ is equal to the trace of  the Hessian of $f$ at $N(z)$.
 \end{proposition}
 
In particular, from the proposition follows that if $f$ is a harmonic function then $H_f$ is a minimal surface. It is interesting to notice that the minimality of $H_f$ follows also from the results of H. Lewy,
see the paper~\cite{Lewy1968}.
 
\begin{proposition}\label{prop2} Let $v$ be a harmonic function defined in a domain $G\ss \R^3$. Suppose
  that $\rank (Hess\, v)=2$ in $G$. Then the set $\nabla v(G)$ is a minimal surface in $\R^3$.
\end{proposition}
   
Proposition \ref{prop2}  also is a consequence of a deep theory of special Lagrangian manifolds, see
 the paper~\cite{HarveyLawson1982}. 
 
 It is interesting to notice that by a remarkable observation of 
M.~Trai\-zet~\cite{Traizet2014}
 the entire solutions of the one-phase free boundary problem \eqref{0} on the plane are related to complete minimal surfaces in $\R^3$. Traizet constructed a Weierstrass-type map from entire solutions of \eqref{0} to immersed minimal surfaces in $\R^3$. Surfaces constructed by Traizet are symmetric with respect to a plane in $\R^3$ and hence they meet that plane orthogonally. It appears that for the simply connected $G$ there are only two entire solutions of \eqref{0} with the
 corresponding minimal surfaces a plane and a catenoid.
 
Notice that a restriction of homogeneous order $1$ harmonic function $u$  on the sphere $\Sp^2$ is an eigenfunction of the Laplacian  on $K\cap\Sp^2$ with the eigenvalue $2$. 
 Thus we can set Theorem \ref{th2} with some generalizations in terms of overdetermined spectral problem. Let $\Om$ be a bounded two dimensional simply connected Riemannian surface of a constant Gaussian curvature $\kappa$ and with a smooth boundary $\partial \Om$. For the Laplace-Beltrami operator $\De$
 on $\Om$ suppose $v$ be a solution of the following overdetermined spectral problem:
 \begin{equation}\label{3} \left\{
\begin{array}{l l}
\De v=  \la v&\mbox{in $\Om$}, \\
v=\al &\mbox{on  $\pa \Om$}, \\
|\nabla v|=\be  &\mbox{on  $\pa \Om$},\\
\end{array} \right. \end{equation}

It is expected that nontrivial solutions of system~\eqref{3} exist only in a disk. 
In the flat case and $\beta =0$  that conjecture is known as the Schiffer's conjecture. Its generalization (generalized Schiffer's conjecture) was widely discussed, see~\cite{Schaefer2001}.  It has a dual integral-geometrical setting, \cite{WillmsChamberlandGladwell1995}.
In the plane case ($\kappa=0$) for $\beta=0$ the above conjecture is equivalent to a long standing  Pompeiu conjecture, see details in a beautiful survey of L. Zalcman~\cite{Zalcman1980}.  On symmetric spaces
the problem was discussed in the paper~\cite{BerensteinZalcman1980}. The case of the unbounded domain $\Om$
was discussed in the paper~\cite{BerestyckCaffarelliNirenberg1997}. For a flat
unbounded $\Om$ and $\lambda=\alpha=0$ a nontrivial example of a solution of \eqref{3} comes from the catenoid
via the Traizet map. However, for bounded solutions $v$ the conjecture holds,
see the paper~\cite{RosRuizSicbaldi2017}.

We prove the following 
   \begin{theorem}\label{th6} Assume $\la= -2\kappa $ and $v$ is a solution of 
   system~\eqref{3}.Then $\Om$ is a geodesic disk. \end{theorem}

\section{Proofs of the theorems}

\begin{lemma}\label{K-constant}
If the Gaussian curvature of a free-boundary mininal surface in a three-dimensional ball
is constant on a connected component of a boundary, then this
component of the boundary is an arc of a circle.
\end{lemma}
 
\noindent {\bf Proof.} We can assume that the ball is of radius $1.$ 
Let $N$ denote a unit normal field on the surface.
Let us choose a point $p$ on the component of the boundary. 
Since the surface meets the boundary sphere $\mathbb{S}^2=\partial\mathbb{B}_1^3$
orthogonally, one can choose such an orthonormal basis in the
three-dimensional space that (i) $N|_p=(0,0,1),$ (ii) the unit tanget vector to the component of
the boundary at $p$ is $(0,1,0),$ (iii) the outward normal vector to $\mathbb{S}^2$ at $p$
is $(1,0,0).$ Remark now that if we put the origin at the center of the ball then
$p=(1,0,0)$ and the sphere is just the standard unit sphere centered at the origin.

Let us parametrize the surface as $(x,y,f(x,y)).$ Then one has
\begin{equation}\label{f-derivatives}
f(1,0)=0,\quad f_x(1,0)=f_y(1,0)=0.
\end{equation}
Then the component of the boundary can be parametrized as
\begin{equation}\label{omega}
\omega(t)=(g(t),t,f(g(t),t)).
\end{equation}
Since $\omega'(0)=(0,1,0)$ and
\begin{equation}\label{omega-derivative}
\omega'(t)=(g'(t),1,f_x(g(t),t)g'(t)+f_y(g(t),t)),
\end{equation}
one has
\begin{equation}\label{g-derivative}
g(0)=1,\quad g'(0)=0.
\end{equation}
Since the surface meets the sphere orthogonally, at each point $\omega(t)$
the unit normal vector
$$
N=\frac{(-f_x,-f_y,1)}{\sqrt{1+f_x^2+f_y^2}}
$$ is
orthogonal to the radius vector of this point, i.e.
\begin{equation}\label{N-radius}
-f_x(g(t),t)g(t)-f_y(g(t),t)t+f(g(t),t)=0.
\end{equation}
If one takes the derivative of equation ~\eqref{N-radius} with respect to $t$ and
one substitutes $t=0,$ then one obtains $f_{xy}(1,0)=0$ due to formulae~\eqref{g-derivative}.

This implies that
$$
H|p=\frac{f_{xx}(1,0)+f_{yy}(1,0)}{2},\quad K|p=f_{xx}(1,0)f_{yy}(1,0).
$$
Since the surface is minimal, one has $H=0.$ It follows that
\begin{equation}\label{fyy}
f_{yy}(1,0)=\sqrt{-K}.
\end{equation}

Let us take
$$
l=\int\limits_0^t|\omega'(t)|dt=\int\limits_0^t\sqrt{(g'(t)^2+1+(f_x(g(t),t)g'(t)+f_y(g(t),t))^2}\,dt
$$
as a natural parameter on the component of the boundary. Then equations~\eqref{f-derivatives}
and~\eqref{g-derivative} imply
$$
\frac{dl}{dt}(0)=1,\quad \frac{d^2l}{dt^2}(0)=0.
$$
It follows that
\begin{equation}\label{d-d}
\frac{d^2\omega}{dt^2}(0)=\frac{d^2\omega}{dl^2}(0)\left(\frac{dl}{dt}(0)\right)^2+
\frac{d\omega}{dl}(0)\frac{d^2l}{dt^2}(0)=\frac{d^2\omega}{dl^2}(0).
\end{equation}

Let us take the derivative of equation~\eqref{omega-derivative}, 
substitute $t=0$ and use \eqref{f-derivatives} and \eqref{g-derivative}. One
obtains
\begin{equation}\label{omega-2d}
\frac{d^2\omega}{dt^2}(0)=(g''(0),0,f_{yy}(1,0)).
\end{equation}
Since the boundary belongs to the sphere, equation~\eqref{f-derivatives} imples 
$$
g(t)^2+t^2+f(g(t),t)^2=1.
$$
Let us take the second derivative and substitute $t=0.$ One obtains 
\begin{equation}\label{g-2}
g''(0)=-1.
\end{equation}

Let us now compute the curvature of the boundary at the point $p.$
Equations~\eqref{d-d}, \eqref{omega-2d}, \eqref{g-2} and \eqref{fyy} imply
$$
k|_p=\left|\frac{d^2\omega}{dl^2}(0)\right|=\left|\frac{d^2\omega}{dt^2}(0)\right|=
\sqrt{(g''(0))^2+(f_{yy}(1,0))^2}=\sqrt{1-K}.
$$

Since $K$ is a constant on the connected component of the boundary, the
curvature of this component is also a constant. It is well known that a curve
of constant curvature on a sphere is an arc of a circle. This finishes the proof.
$\Box$

\noindent {\bf Proof of Theorem~\ref{th1}.}
Let $A$ be the annulus~\eqref{annulus}, and
$\langle\cdot\,,\cdot\rangle$ be the standard $\mathbb{C}$-bilinear scalar product on $\mathbb{C}^3.$

Consider a map $u: \bar{A}\longrightarrow\mathbb{\bar{B}}^3_1\subset\mathbb{R}^3$
such that $u\in C^2(A)\cap C^1(\bar{A}),$ $u(\partial A)\subset\mathbb{S}^2=\partial\mathbb{B}_1^3,$
\begin{equation}\label{minimal}
u_{z\bar{z}}=0,\quad\langle u_z,u_z\rangle=0,
\end{equation}
and $u(A)$ meets $\mathbb{S}^2$ orthogonally.

Denote by $N$ the unit normal field on $u(A)$ and by
$u^\perp_{zz}$ the component of $u_{zz}$ normal to $u(A)$, i.e.
$$
u_{zz}=u^\perp_{zz} +\frac{\langle u_{zz},u_{\bar{z}}\rangle}{\langle u_z,u_{\bar{z}}\rangle}u_z+
\frac{\langle u_{zz},u_z\rangle}{\langle u_z,u_{\bar{z}}\rangle}u_{\bar{z}}.
$$
In fact, the second formula from~\eqref{minimal} implies that $\langle u_{zz},u_z\rangle=0.$

Consider polar coordinates $r,\theta$ such that $z=re^{i\theta}.$ Since $u_r$ and
$u_\theta$ are tangent to $u(A),$ one has
\begin{equation}\label{uzz}
u^\perp_{zz}=\frac{e^{-2i\theta}}{4}\left(u^\perp_{rr}-\frac{2i}{r}u^\perp_{r\theta}-
\frac{1}{r^2}u^\perp_{\theta\theta}\right).
\end{equation}
The free boundary condition, i.e. the condition that $u(A)$ meets $\mathbb{S}^2$ orthogonally,
implies $u_{r}|_{\partial A}\perp\mathbb{S}^2.$ Hence, one has
$u_r=\lambda u$ on $\partial A$ for some function $\lambda.$ It follows that
$$
u_{r\theta}=\lambda_\theta u+\lambda u_\theta=\frac{\lambda_\theta}{\lambda}u_r+\lambda u_\theta
$$
is a tangent vector. This means that $u^\perp_{r\theta}=0$ on $\partial A.$
Then equation~\eqref{uzz} implies that
$e^{4i\theta}\langle u^\perp_{zz},u^\perp_{zz}\rangle$ is real on $\partial A.$
Since $z^4=r^4e^{4i\theta},$ it follows that $z^4\langle u^\perp_{zz},u^\perp_{zz}\rangle$ is real
and positive on $\partial A.$

It is well-known that in a simple connected domain for a minimal surface $u$ there exists
an adjoint surface $u^*$ such that $f=u+iu^*$ is holomorphic and 
$l=\langle f'',N\rangle$ is a holomorphic function (including the branch points), 
see e.g.~\cite{DierkesHildebrandtSauvigny2010}. It follows that
\begin{equation}\label{z4uzz}
z^4\langle u^\perp_{zz},u^\perp_{zz}\rangle=\frac{z^4}{4}l^2\end{equation}
is also a holomorphic function. $A$ is not simply connected, but the
property of being holomorphic is local, one can check it in simple
connected neighbourhoods of points of $A.$ Hence, $z^4\langle u^\perp_{zz},u^\perp_{zz}\rangle$
is holomorphic on $A.$ Since $z^4\langle u^\perp_{zz},u^\perp_{zz}\rangle$ is real
on $\partial A,$ this function is constant on $A,$
\begin{equation}\label{z4uzz-l}
z^4\langle u^\perp_{zz},u^\perp_{zz}\rangle=k\in\mathbb{R}, \quad k>0.
\end{equation}

Let us now consider the celebrated Enneper-Weierstrass representation.  For $\Phi=2u_z$
there exist a holomorphic function $\mu$ and meromorphic function $\nu\not\equiv0$ such that
$\mu\nu^2$ is holomorphic and
$$
\Phi=\left(\frac{1}{2}\mu(1-\nu^2),\frac{i}{2}\mu(1+\nu^2),\mu\nu\right).
$$
It follows that
$$
u(z)=u_0+\Real\int_{z_0}^z\Phi(z)\,dz.
$$
Let $z=x+iy.$
Let us recall that
\begin{equation}\label{ds2}
ds^2=\Lambda(dx^2+dy^2), \quad \Lambda=\frac{1}{4}|\mu|^2(1+|\nu|^2)^2,
\end{equation}
\begin{equation}\label{l-mu-nu}
l=\langle f'',N\rangle=-\mu\nu',
\end{equation}
see e.g.~\cite{DierkesHildebrandtSauvigny2010}. Remark that $f$ could be multivalued since
$A$ is not simply connected, but $f'=\Phi$ is a holomorphic function.

Let us recall that $z_0$ is a branch point if and only if $z_0$ is a zero of $\mu$ and $\mu\nu^2.$
Since there is no branch points on the boundary, see e.g.~\cite{ColdingMinicozzi2011},
$\mu$ has no zeroes on $\partial A.$

Let us now consider the point $z_0=1$ or $z_0=\rho$ on $\partial A.$ Remark that any point on
$\partial A$ could be transformed by a rotation to one of these two points.
Consider the curve $u(|z_0|e^{i\theta})\subset\mathbb{S}^2$ parametrized by $\theta.$

It is well known that for a curve lying on a sphere its osculating spheres coincide with the initial sphere.
It is also well known that the circle obtained as intersection of the osculating sphere and the osculating
plane at a point of a curve touches this curve at the second order at this point. It follows that there
exist a circle $\gamma(t)\subset\mathbb{S}^2,$  where $t$ is an affine natural parameter, such that
$$
u(z_0)=\gamma(0),\quad u_\theta(z_0)=\dot{\gamma}(0),\quad u_{\theta\theta}(z_0)=\ddot{\gamma}(0).
$$
Performing, if necessary, a rotation and reflection of $\mathbb{R}^3,$ we can assume that
(i) the circle $\gamma$ is the circle $\gamma^3(t)=const,$ (ii) $\gamma^2(0)=0$, here the superscripts
mean the coordinate number, (iii) $N(z_0)\ne(0,0,1).$
Let us remark that property (i) implies 
\begin{equation}\label{u3}
u^3_\theta(z_0)=0,\quad u^3_{\theta\theta}(z_0)=0,
\end{equation}
property (ii) implies
\begin{equation}\label{u1}
u^1_\theta(z_0)=0,
\end{equation}
and (iii) implies that $\nu$ does not have a pole at $z_0,$ see e.g.~\cite{DierkesHildebrandtSauvigny2010}.

Combining
\begin{equation}\label{u_theta}
u_\theta=izu_z-i\bar{z}u_{\bar{z}}=-\Imaginary(z\Phi)
\end{equation}
with equations~\eqref{u3} and~\eqref{u1}, one has
$$
u_\theta(z_0)=\left(0,-\Real\left(\frac{z_0}{2}\mu(1+\nu^2)|_{z=z_0}\right),0\right).
$$
Computing $|u_\theta(z_0)|$ directly and using formula~\eqref{ds2}, we obtain the equation
$$
\left|\Real\left(\frac{z_0}{2}\mu(z_0)(1+\nu^2(z_0))\right)\right|=
\frac{|z_0|}{2}|\mu(z_0)|(1+|\nu(z_0)|^2)|.
$$
It is easy to prove that for $a,b\in\mathbb{C}$ the inequality 
$$
|\Real(a(1+b^2))|\leqslant|a|(1+|b|^2)
$$ holds, and one has the equality
if and only if $a,b\in\mathbb{R}$ or $a=0.$

Since $\mu$ has no zeroes on $\partial A,$
it follows that $\mu(z_0),$ $\nu(z_0)\in\mathbb{R}.$

Since $\gamma(t)$ is a circle parametrized by an affine natural parameter,
$(\ddot{\gamma}^1(0),\ddot{\gamma}^2(0))\parallel (\gamma^1(0),\gamma^2(0)).$
Then $\gamma^2(0)=0$ implies 
$$
u^2_{\theta\theta}(z_0)=\ddot{\gamma}^2(0)=0.
$$
Since
$$
u_{\theta\theta}=-\Real(z\Phi+z^2\Phi'),
$$
one has
$$
u^2_{\theta\theta}(z_0)=-\Real\left(z\frac{i}{2}\mu(1+\nu^2)+z^2\frac{i}{2}(\mu'(1+\nu^2)+2\mu\nu\nu')\right)|_{z=z_0}=0.
$$
It follows that
$$
(\mu'(1+\nu^2)+2\mu\nu\nu')|_{z=z_0}\in\mathbb{R}.
$$
But formula~\eqref{z4uzz-l} implies $\mu\nu'|_{z=z_0}\in\mathbb{R}.$
It follows that $\mu'(z_0),$ $\nu'(z_0)\in\mathbb{R}.$

Let us remark that
$$
2\mu(z_0)|\mu|_\theta(z_0)=\frac{\partial}{\partial\theta}|\mu|^2(z_0)=
(\mu_\theta\bar{\mu}+\mu\bar{\mu}_\theta)|_{z=z_0}=
(2\Real\mu_\theta(z_0))\mu(z_0).
$$
Since $\mu_\theta(z_0)=iz\mu'(z_0)$ is purely imaginary and $\mu(z_0)\ne0,$
one has $|\mu|_\theta(z_0)=0.$ The same argument proves that at least one
of quantities $|\nu|_\theta(z_0)$ or $\nu(z_0)$ is zero.

Consider now $\sqrt\Lambda=\frac{|\mu|}{2}(1+|\nu|^2).$ One has
$$
\frac{\partial}{\partial\theta}\sqrt\Lambda(z_0)=
\left(\frac{|\mu|_\theta}{2}(1+|\nu|^2)+|\mu||\nu||\nu|_\theta\right)|_{z=z_0}=0.
$$
Since the formula $\Lambda_\theta=0$ does not change under rotations of $z$-plane or
isometries of $\mathbb{R}^3,$ it holds for any point on $\partial A.$
This means that $\theta$ is an affine natural parameter on each
connected component of $u(\partial A).$

Moreover, the formula for Gaussian curvature, see e.g.~\cite{DierkesHildebrandtSauvigny2010},
$$
K=-\left(\frac{4|\nu'|}{|\mu|(1+|\nu|^2)^2}\right)^2
$$
and equations~\eqref{z4uzz}, \eqref{ds2} and~\eqref{l-mu-nu}  imply that
$K=-\frac{4k}{|z|^4\Lambda^2}$ is a constant
on each
connected component of $u(\partial A).$

Then Lemma~\ref{K-constant} implies that each component of $u(\partial A)$
is a circle lying on $\mathbb{S}^2$. 

Let a circle $\sigma$ be a boundary component of $u(A)$ and  $p\in\sigma$. Let us consider the case
when $\sigma$ is not a great circle on $\mathbb{S}^2.$ In this case one can find a catenoid $H$ which meets
orthogonally $S^2$ at $\sigma$. Without loss of generality we can assume that the initial minimal surface
could be locally parametrised as $(x,y,f(x,y))$ and the catenoid as $(x,y,h(x,y)).$
Since $f$ satisfies the minimal surface equation and have the same Cauchy
data on $\sigma$ in a neighbourhood of $p$ as $h,$ one has $f(x,y)=h(x,y)$ where 
both functions are defined. Then $u(A)$ and $H$ coincide globally and $u$
is a reparametrisation of a catenoid. Since all conformal automorphisms
of the annulus $A$ are described by Schottky theorem~\cite{Schottky1877},
$$
z \mapsto \lambda z^{\pm 1},
$$
where $|\lambda|=1$ or $\rho,$ we obtain that $u$ is a catenoid. 

In the case when $\sigma$ is a great circle the same argument could be applied with
a plane disk instead of a cathenoid. As a result, another connected component
of $u(\partial A)$ is inside the ball. This contradicts the assumption $u(\partial A)\subset\mathbb{S}^2.$
This finishes the proof.
$\Box$

Assuming a non-zero Dirichlet boundary condition on $\pa K$ we will consider a generalization of the problem \eqref{1}: 

\begin{equation}\label{pr2} \left\{
\begin{array}{l l}
\De v=  0&\mbox{in $K$}, \\
v=\al |x|&\mbox{on  $\pa K$}, \\
|\nabla v|=1  &\mbox{on  $\pa K\sm \{ 0\}$},\\
\end{array} \right. \end{equation}
where $v$ is a homogeneous degree $1$ function. We assume that $\al\in \R$ is a constant. Then the second boundary condition implies that $\al \in (-1,1)$.
 
 Let $v$ be a homogeneous order $1$ harmonic function defined in the cone $K$ and satisfying equation \eqref{pr2}. Denote $G=K\cap \Sp^2,\, f=\nabla v: G\to \R^3$. Let $\ga\ss G$ be the set of critical points of $f$, i.e., the set of vanishing of the differential $df$. Since $f$ is a real analytic map, $\ga$ is either a set of isolated points in $G$, or it contains a one dimensional ark $\ga'$, see the paper~\cite{Whitney1965}. Consider  the second case. Let $\tilde \ga$ be the conic extension of $\ga'$ to $\R^3$.Then 
there is a linear function $l$ in $\R^3$ such that $l-v=0, \, \nabla v-\nabla l=0$ on $\tilde \ga$. Since $l-v$ is a harmonic function then by the uniqueness of the
solution of the Cauchy problem $l=v$ in $K$ and the theorem follows. 
Suppose now that $\ga'$ is a set of isolated points.
 By Proposition \ref{prop} surface $\nabla v (G\sm \ga') \ss \R^3$ is a minimal surface. By a theorem by 
 Gulliver and Lawson, see~\cite{GulliverLawson1986} and~\cite{Meeks2007}, the surface $\Ga$ could be
 extended to the set $\ga'$ as a branching minimal surface $\Ga= \nabla v( G)$. Assume now that $v$ satisfies in $K$ equation \eqref{1}. Since $|\nabla v|=1, \, \pa \Ga \ss \Sp^2$. Let $x\in \pa K$. Since $v$ is a homogeneous order $1$ function
 $\nabla v (x) =\al x +\beta e$, where $|e|=1, (x,e)=0$, $\sqrt{\al^2+\beta^2}=1$. Thus the angle between vectors $x$ and $\nabla v(x)$ is fixed for all points $x\in \pa K$.
 By Proposition \ref{prop} $x$ is a unit normal to $\Ga$ at $\nabla v(x)$
 and hence the angle between vectors $N(\nabla v(x))$ and $\nabla v(x)$ is fixed for all points $x\in \pa K$. Thus $\Ga$ intersects  sphere $\Sp^2$ under a fixed angle in particular if $\al=0$ then $\Ga$ meats  sphere $\Sp^2$
 orthogonally.  Now Theorem \ref{th2}  follow from the theorem of Nitsche and from Theorem \ref{th1}.
 
 {\bf Remark.} It is easy to see that the curves $\partial \Ga$ and $\partial K\cap \Sp^2_r$ are dual curves on
 the sphere $\Sp^2_r$.
 
{\bf Proof of Theorem \ref{th6}.}  Since the Gaussian curvature of $\Om$
 is $1$ there exists an isometry
 $$
 i: \Om \longrightarrow \Sp^2
 $$ 
and  $i(\Om)$ is a domain possibly multi-sheeted on $\Sp^2$. Denote by $u$ the pull down of the function $v$
 from $\Om$ to $i(\Om)\ss\Sp^2$. We will assume that the function $u$ is extended as a homogeneous $1$ function 
to a cone $K\ss \R^3$ over $i(\Om)$, where $K$ is possibly multi-sheeted cone. Then $u$ satisfies equation \eqref{pr2}.  
The same argument as above shows that the surface $\nabla u$ intersects  sphere $\Sp^2_r$ under a fixed angle.
Hence Theorem \ref{th6} follows from Nitsche's theorem.

\end{document}